\documentclass[12pt]{article}
\newcommand {\beq} {\begin{equation}}
\newcommand {\eeq} {\end{equation}}
\newcommand {\beqa} {\begin{eqnarray*}}
\newcommand {\eeqa} {\end{eqnarray*}}
\newcommand {\bc} {\begin{center}}
\newcommand {\ec} {\end{center}}

\newcommand {\VE} {\varepsilon}
\newcommand {\Om} {\Omega}
\newcommand {\pOm} {\partial\Omega}
\newcommand {\vb} {\vec\beta}

\newcommand {\rmon} {\ {\rm on} \ }

\begin{document}
\title{ Numerical methods for convection-diffusion problems\\[.5em]
or\\[.5em]
The 30 years war}

\author{Martin Stynes\footnote{This is the text of an invited lecture from the 20th Biennial Numerical Analysis held in Dundee, Scotland during 24-27 June 2003. The Proceedings of this conference were published electronically on the website of the University of Dundee in 2003, but seem to be no longer available. As I get occasional requests for copies of this paper, I am now placing this copy of it on the arXiv server.}\\Department of Mathematics, National University of Ireland,
Cork, Ireland\\ m.stynes@ucc.ie}

\maketitle

\begin{abstract}
Convection-diffusion problems arise in the modelling of many
physical processes. Their typical solutions exhibit boundary
and/or interior layers. Despite the linear nature of the
differential operator, these problems pose still-unanswered
questions to the numerical analyst.

This talk will give a selective overview of numerical methods for
the solution of convection-diffusion problems, while placing them
in a historical context. It examines the principles that underpin
the competing numerical techniques in this area and presents some
recent developments.

\end{abstract}

\thispagestyle{empty}

\section{Talk overview}

To quote the opening words of Morton's book \cite{Mo96}: \lq\lq
Accurate modelling of the interaction between convective and
diffusive processes is the most ubiquitous and challenging task in
the numerical approximation of partial differential equations."
% I trust this is sufficient justification for my talk!

I shall describe the nature of (steady-state) convection-diffusion
problems, then draw some comparisons between the development of
numerical methods for convection-diffusion problems during the
last 30 years and that well-known 17th-century conflict known as
the 30 Years War, whose history most Europeans learn during their
schooldays.

Let's begin by reminding ourselves of its main features.

\section{The phases of the 30 Years War}
The 30 Years War war began in 1618 as a struggle between certain
Catholic and Protestant states, but eventually sucked in all the
major European countries. It devastated much of present-day
Germany and the Czech Republic; indeed, not until the 2nd World
War, 300 years later, was there a comparable amount of destruction
in any European war.

As various protagonists entered or left the the conflict, the
action passed through several {\em phases}
\cite[pp.252--255]{Peng1}: the 1618--23 Bohemian phase, the
1624--29 Danish phase, the 1630--35 Swedish phase, and the
1635--48 Franco-Swedish phase.

To indicate how science was progressing at this time, we note that
Kepler's 3rd law (proportionality of the square of the period of
revolution of a planet to the cube of the length of the major axis
of its orbit) was published in 1619. He was then living in Linz,
in Austria.

The war ended with the Peace of Westphalia (Westfalen in German)
in 1648.

\section{Convection-diffusion problems}\label{cdprobs}
These take the form
\begin{equation}\label{cd}
 -\VE \Delta u +  \vb\cdot\nabla u  = f \rmon \Om,
\end{equation}
with some boundary conditions on $\pOm$, where $-\VE\Delta u$ models
diffusion and $\vb\cdot\nabla u$ models convection. The parameter
$\VE$ is positive but small; think of it as say $10^{-6}$. The
region $\Omega$ is any reasonable domain in $n$ dimensions, where
$n\ge 1$. The differential operator is elliptic, so under suitable
hypotheses on $\vb$ and the boundary data, (\ref{cd}) has a
solution in $C^2(\Omega)$. Here we take $\vb \approx O(1)$, i.e.,
{\em convection dominates diffusion}:
\[
\frac{|{\rm coefficient\ of\ }\nabla u|}{|{\rm coefficient\ of\ }\Delta u|}
= \frac{|\vb|}{|\VE|} >\! > 1.
\]

For most boundary conditions, this is an example of a singularly
perturbed partial differential equation (PDE).
Convection-diffusion PDEs arise in many applications \cite{Mo96}
such as the linearized Navier-Stokes equations and the
drift-diffusion equation of semiconductor device modelling.

To get some feeling for the behaviour of solutions to such
problems, let's consider a simple example in one dimension, that
is, an ordinary differential equation (ODE):
  \beqa
-\VE u''  + u' &=& 2 \quad {\rm on\ } (0,1),\\
u(0) = u(1) &=& 0.
 \eeqa
Then
\begin{eqnarray}
u(x) & = & 2x + \frac{2(e^{-1/\VE} - e^{-(1-x)/\VE})} {1 -
e^{-1/\VE}} \nonumber\\
 & = & 2x - 2e^{-(1-x)/\VE} + O( e^{-1/\VE}). \label{decomp}
\end{eqnarray}
Here $2x$ is the solution of the first-order problem $u'(x) = 2, \
u(0) = 0$; the rapidly-decaying exponential $e^{-(1-x)/\VE}$ is a
boundary layer function---i.e., it is not large but its
first-order derivative is large near $x=1$; and the final term
$O(e^{-1/\VE})$ is negligible.

For PDEs also, the solution $u$ of (\ref{cd}) has a structure
analogous to (\ref{decomp}): it can be written as the sum of the
solution to a 1st-order hyperbolic problem + layer(s) + negligible
terms. Let's make this a little more precise. Divide the boundary
$\pOm$ into 3 parts:
 \beqa {\rm inflow\ boundary\ }
\partial^-\Omega & = & \{
  x\in\partial\Omega: \vec\beta.\vec n < 0 \}, \\
  {\rm outflow\ boundary\ }\partial^+\Omega
    & = & \{ x\in\partial\Omega: \vec\beta.\vec n > 0 \}, \\
{\rm tangential\ flow\ boundary\ }\partial^0\Omega & = & \{
x\in\partial\Omega: \vec\beta.\vec n = 0 \}, \eeqa
where $\vec n$ is the outward-pointing unit normal to $\pOm$. Then
the 1st-order hyperbolic problem is $\vb\cdot\nabla u = f$ on
$\Om$, with boundary data specified on $\partial^-\Omega$.
Professor Ron Mitchell, one the founders of the Dundee
conferences, used to refer to the difficulties one can face in
attempting to solve accurately this innocent-looking problem as
\lq\lq the great embarrassment of numerical analysis." Usually
(depending on the precise boundary conditions in (\ref{cd})) the
solution $u$ has an {\em exponential boundary layer} along
$\partial^+\Omega$ and {\em parabolic/characteristic boundary
layers} along $\partial^0\Omega$. Exponential layers are, at each
point in $\partial^+\Omega$, essentially the same as the function
$e^{-(1-x)/\VE}$, but characteristic layers have a much more
complicated structure and cannot be defined by an ODE. They
nevertheless have the layer quality of fast decay in a narrow
region (roughly of width $O(\sqrt\VE\,)$ along $\partial^0\Omega$.

There will also be a characteristic layer in the interior of $\Om$
emanating from each point of discontinuity in the boundary
conditions on $\partial^- \Omega$ (think how a discontinuity would
be propagated across $\Om$ by $\vb\cdot\nabla u = f$; the effect
of the diffusion term $-\VE \Delta u$ is to smooth this discontinuity
into a continuous but steep layer).

For some illustrations of how solutions to such problems can look,
see, e.g., \cite{MaSt96, MaSt97}.

\section{Numerical instability}\label{NI}
Standard numerical approximations of differential equations use a
central difference approximation of the convective term. That is,
for ODEs, one approximates $ u'(x_i)$ by $(u_{i+1}^N -
u_{i-1}^N)/(2h)$ in the usual notation, where $h$ is the local
mesh-width. On quasiuniform meshes this yields oscillatory and
inaccurate solutions; see, e.g., \cite[Fig.~5]{HB79}.

One can give many indirect explanations of this poor performance.
For example, a careful inspection of the analysis of standard
numerical methods reveals an assumption that the diffusion
coefficient is bounded away from 0; but in (\ref{cd}) the
parameter $\VE$ can be very small, and this means that the
standard analysis is no longer valid. An alternative explanation:
the matrices generated by the approximations of the differential
operator are not $M$-matrices when $\VE$ is small relative to the
local mesh-width, so their inverses can be expected to have both
positive and negative entries. As a consequence the computed
solution will display oscillations.

One means of eliminating oscillations is to approximate the
convective derivative by a non-centered approximation called {\em
upwinding}: when solving $-\VE u'' + u' = f$ on a uniform mesh of
width $h$, replace
\[
u'(x_i)  \mapsto  (u_{i+1}^N - u_{i-1}^N)/(2h)
\]
by
\[
u'(x_i)  \mapsto  (u_{i}^N - u_{i-1}^N)/h,
\]
while discretizing $-\VE u''(x_i)$ in the usual way. (Here
$\{u_i^N\}_{i=0}^N$ is the computed solution.) That is, the
approximation of $u'(x_i)$ uses values of $u_i^N$ that are chosen
away from the boundary layer at $x=1$. It is easy to check that
this modified discretization of $-\VE u'' + u'$ yields an
$M$-matrix; consequently the computed solution is more stable and
no longer has non-physical oscillations.

While the unwanted oscillations have disappeared, this has come at
a price: layers in the computed solution are excessively smeared,
i.e., are not as steep as they should be. See, e.g.,
\cite[Fig.~5]{HB79}. To motivate a way of addressing this
shortcoming, we observe that on a uniform mesh of width $h$,
upwinding yields
 \beqa
 (-\VE u''+u')(x_i) & \mapsto  &
\frac{-\VE}{h^2}(u_{i+1}^N -2 u_i^{N} + u_{i-1}^N)
 + \frac{1}{h} (u_i^N - u_{i-1}^N) \\
& = &\hspace{-0.1cm}- \left( \VE + \frac{h}{2}\right)
\frac{1}{h^2}(u_{i+1}^N -2 u_i^{N} + u_{i-1}^N) + \frac{1}{2h}
(u_{i+1}^N - u_{i-1}^N).
 \eeqa
That is, upwinding applied to $-\VE u'' + u'$ is the same method
as standard central differencing applied to $-(\VE + h/2) u'' +
u'$. To put this in words, we can regard upwinding as the standard
discretization of a modified differential equation---modified by
artificially increasing the diffusion coefficient by $h/2$.

Now we see the possibility of modifying the diffusion coefficient
by some other quantity before applying a standard numerical
method, with the aim of retaining stability while introducing less
smearing of layers in the computed solution. This way of thinking
turns out to be quite fruitful; in fact, stable numerical methods
on uniform meshes for convection-diffusion ODEs are usually
equivalent to modifying the diffusion in the original differential
equation then applying a standard method (e.g., central
differencing)---but for PDEs, the connection may be less
straightforward.

{\em Summary:} when a standard numerical method is applied to a
convection-diffusion problem, if there is too little diffusion,
then the computed solution is oscillatory, while if there is too
much diffusion, the computed layers are smeared.

One can add artificial diffusion using finite difference, finite
element or finite volume methods. See \cite{RST96} for many
illustrations of how this can be done. In the rest of this talk, I
shall discuss a few well-known techniques for the numerical
solution of convection-diffusion problems  that operate in this
way.

\section{1969--early 1990s: the international phase}\label{intl}
Our history of numerical methods for convection-diffusion problems
begins about 30 years ago, in 1969. In this year, two significant
Russian papers \cite{Ba69, Il69} analysed new numerical methods
for convection-diffusion ODEs.

In \cite{Ba69}, Bakhvalov considered an upwinded difference scheme
on a layer-adapted graded mesh. Such meshes are based on a
logarithmic scale (the inverse of the exponential layer function
that we met in (\ref{decomp})). They are very fine inside the
boundary layer and coarse outside. The fineness of the mesh means
that the added artificial diffusion is very small inside the
layer, and consequently the layer is not smeared excessively.

We shall return later to Bakhvalov's idea, as initially it was
less influential than \cite{Il69}, where A.M.Il'in  used a uniform
mesh but chose the amount of added artificial diffusion in such a
way that for constant-coefficient ODEs the computed solution
agrees exactly with the true solution at the meshpoints. The
amount of artificial diffusion involves exponentials, and schemes
of this type are called {\em exponentially-fitted} difference
schemes. See \cite{Ro} or \cite{RST96} for details of the scheme.
(In fact the same scheme had been used much earlier in \cite{AS55}
but no analysis of its behaviour was given there.)

During the next 20 years, researchers from many countries
developed Il'in-type schemes for many singularly perturbed ODEs
and some PDEs. See \cite{RST96} for references; here we just
mention Griffiths and Mitchell from Dundee.

The original Il'in paper used a complicated technique called the
\lq\lq double-mesh principle" to analyse the difference scheme.
This became obsolete overnight when in 1978 Kellogg and Tsan
published a revolutionary and famous paper \cite{KT78} that was
gratefully seized on by other researchers in the area. Their paper
showed how to design {\em barrier} or {\em comparison functions}
to convert truncation errors to computed errors, and also gave for
the first time sharp a priori estimates for the solution of the
convection-diffusion ODE. (Historical note: it was the first of
many papers by Bruce Kellogg on convection-diffusion problems, and
it was the only mathematical paper that Alice Tsan ever wrote!)

Today exponential fitting is still used for instance in the
well-known package PLTMG and in semiconductor device modelling
(where it's known as the Scharfetter-Gummel scheme).  In
\cite{An00}, Angermann gives an example of an exponentially-fitted
scheme that does a remarkable job of capturing an interior layer
on a uniform mesh.

A related idea is the residual-free bubbles FEM that has been
developed and analysed in recent years by Brezzi, Franca, et al.
This doesn't explicitly contain exponentials, but it is based on
the idea of solving a local problem exactly \cite{FNS98}, as is
Il'in's method \cite{Ro}.

\section{1979--mid 1990s: the Swedish phase}

The work described in the previous section is finite difference in
nature. In 1979, Hughes and Brooks \cite{HB79} introduced the {\em
streamline diffusion finite element method} (SDFEM) for
convection-diffusion problems. This kick-started a development of
finite element methods for convection-diffusion problems that
continues to this day. Many researchers have participated in this
effort, but from the point of view of analysis, the dominant
character has been Claes Johnson from Sweden. See \cite{RST96} for
an overview of the relevant literature.

The SDFEM (also known as SUPG) works as follows. Consider the PDE
(\ref{cd}) with homogeneous Dirichlet boundary conditions, and
$\Om$ a convex bounded subset of $R^2$. Assume that $\vb$ is
constant with $|\vb| = 1$. Write for convenience $\vb\cdot\nabla u
\equiv u_\beta$.

Suppose we have a triangular mesh on $\Om$. We'll discuss only
piecewise linears here, but there is an analogous, slightly more
complicated method for piecewise polynomials of higher degree (see
\cite{RST96}). For the piecewise linear trial space $V^N \subset
H_0^1(\Om)$, the SDFEM is: find $u^N \in V^N$ such that
\[
\VE(\nabla u^N, \nabla v^N) + ( (u^N)_\beta, v^N + \delta v_\beta^N)
= (f,v^N + \delta v_\beta^N)
\]
for all $v \in V^N$, where $\delta$ is a user-chosen locally
constant parameter with $\delta \ge 0$; typically $\delta  =
O({\rm local\ mesh\ diameter})$

This is roughly equivalent to altering the PDE from $-\VE\Delta u +
u_\beta = f$ to $ -\VE\Delta u - \delta u_{\beta\beta} + u_\beta  =
f$, then applying the standard Galerkin FEM---i.e., \textit{we
have added artificial diffusion only in the streamline/flow
direction}. This stabilizes the SDFEM and can remove the
outflow-layer oscillations one would obtain if the standard
Galerkin FEM were applied directly to (\ref{cd}); moreover, as
diffusion is added only in the direction of flow, one does not
smear characteristic layers.

How exactly should one choose $\delta$? No \lq\lq optimal" formula
is known. Different choices introduce different amounts of
diffusion. In \cite[Figs. 2, 3]{MaSt96} one sees the striking
effect of different choices of $\delta$ on the sharpness of
computed outflow boundary layers.

In practice the SDFEM yields accurate solutions away from layers
and local error estimates of Johnson, Schatz and Wahlbin, and
Niijima, reflect this behaviour. See \cite{RST96}. On subdomains
$\Om_0$ of $\Om$ that lie \lq\lq away from" layers,
\[
|(u-u^N)(x)| \le C h^{11/8} \ln (1/h)\ ||u||_{C^2(\Om_0)}.
\]
This is almost sharp: numerical results  of Zhou imply that in
general $O(h^{3/2})$ is the best possible bound for piecewise
linears.

But the stabilization of the SDFEM has little effect along
characteristic layers. Kopteva \cite{Ko03} shows that one obtains
only $O(\delta)$ pointwise accuracy inside parabolic boundary and
interior layers, and as $\delta = O(h)$ typically, this means one
can at best get first-order convergence inside characteristic
layers, even on special meshes.

The idea of the SDFEM has generated several related FEMs for
convection-diffusion problems: the Galerkin least-squares FEM,
negative-norm stabilization of the FEM, and the currently popular
discontinuous Galerkin FEM.

\section{1990--present: the Russian-Irish phase}
Recall from Section \ref{intl} that Bakhvalov-type graded meshes
can be used to solve convection-diffusion problems. In 1990 the
Russian mathematician Grisha Shishkin showed that instead one
could use a simpler piecewise uniform mesh. This idea has been
enthusiastically propagated throughout the 1990s by a group of
Irish mathematicians: Miller, O'Riordan, Hegarty and Farrell. See
\cite{FaHe00,RST96} and their bibliographies.

The Shishkin mesh is chosen a priori. It is very fine near layers
but coarse otherwise. For example, if the domain $\Om$ is the unit
square and the problem is
\[
-\VE\Delta u + b_1 u_x + b_2 u_y = f, \quad {\rm with\ }b_1 > 0, b_2
> 0,
\]
so one has boundary layers along $x=1$ and $y=1$, then this
tensor-product mesh has transition points (where the mesh switches
from coarse to fine) at $1-\lambda_x$ and $1-\lambda_y$ on the
$x$- and $y$-axes respectively, where $ \lambda_x = (4\VE/b_1)\ln
N$ and $\lambda_y = (4\VE/b_2)\ln N$. Here $N$ is the number of
mesh points in each coordinate direction. The fine and coarse mesh
regions on the coordinate axes each contain $N/2$ mesh intervals.
See Figure 1 for the mesh with $N=8$ (the mesh rectangles have
been bisected into triangles to permit the use of a piecewise
linear FEM).
\begin{figure}\label{Shishkin}
\unitlength0.5pt
\begin{center}
\begin{picture}(241,241)
 \put(  0,  0){\line(0,1){240}}
 \put( 50,  0){\line(0,1){240}}
 \put(100,  0){\line(0,1){240}}
 \put(150,  0){\line(0,1){240}}
 \put(200,  0){\line(0,1){240}}
 \put(210,  0){\line(0,1){240}}
 \put(220,  0){\line(0,1){240}}
 \put(230,  0){\line(0,1){240}}
 \put(240,  0){\line(0,1){240}}
 \put(  0,  0){\line(1,0){240}}
 \put(  0, 40){\line(1,0){240}}
 \put(  0, 80){\line(1,0){240}}
 \put(  0,120){\line(1,0){240}}
 \put(  0,160){\line(1,0){240}}
 \put(  0,180){\line(1,0){240}}
 \put(  0,200){\line(1,0){240}}
 \put(  0,220){\line(1,0){240}}
 \put(  0,240){\line(1,0){240}}
 \multiput(  0, 40)(50,0){4}{\line(5,-4){50}}
 \multiput(  0, 80)(50,0){4}{\line(5,-4){50}}
 \multiput(  0,120)(50,0){4}{\line(5,-4){50}}
 \multiput(  0,160)(50,0){4}{\line(5,-4){50}}
 \multiput(200, 40)(10,0){4}{\line(1,-4){10}}
 \multiput(200, 80)(10,0){4}{\line(1,-4){10}}
 \multiput(200,120)(10,0){4}{\line(1,-4){10}}
 \multiput(200,160)(10,0){4}{\line(1,-4){10}}
 \multiput(  0,180)(50,0){4}{\line(5,-2){50}}
 \multiput(  0,200)(50,0){4}{\line(5,-2){50}}
 \multiput(  0,220)(50,0){4}{\line(5,-2){50}}
 \multiput(  0,240)(50,0){4}{\line(5,-2){50}}
 \multiput(200,180)(10,0){4}{\line(1,-2){10}}
 \multiput(200,200)(10,0){4}{\line(1,-2){10}}
 \multiput(200,220)(10,0){4}{\line(1,-2){10}}
 \multiput(200,240)(10,0){4}{\line(1,-2){10}}
\end{picture}
\caption{Shishkin mesh with $N=8$}
\end{center}
\end{figure}
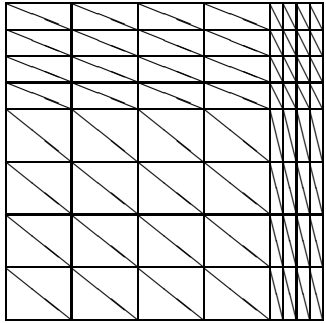
Shishkin and his coauthors favour the use of upwinding (see
Section \ref{NI}) on this mesh. Since the mesh is fine at the
boundary layers, upwinding does not smear these layers. The
computed solution has no non-physical oscillations, and one
usually obtains almost first-order (i.e., up to a factor $\ln N$)
pointwise convergence at the mesh nodes. Unlike the SDFEM, one
does not have to manage a free parameter. The computed solutions
look satisfactory \cite{FaHe00}.

One can of course instead use a FEM on a Shishkin mesh. Lin\ss\
and Stynes \cite{LS01} give numerical results for linears and
bilinears on these meshes, and show that bilinears are more
accurate in the layer regions. Further theoretical evidence that
bilinears are superior to linears is given in \cite{ST03}.

The drawbacks to Shishkin meshes are that one must know the
location and nature of the layers a priori, and up to now the
method has been implemented only on rectangular domains. Curved
interior layers have not been tested numerically using exact
Shishkin meshes, but in \cite{MaSt97} an interior layer is
computed accurately using PLTMG and a simple heuristic
approximation of a Shishkin mesh. Furthermore, the analysis of
schemes on these meshes requires strong assumptions on the data of
the problem to ensure sufficient differentiability of solution $u$
and thereby justify the choice of mesh.

An excellent survey of the published literature for layer-adapted
meshes (Shishkin, Bakhvalov, etc.) applied to convection-diffusion
problems is given by Lin\ss\ \cite{Li03}.

\section{A historical connection}

If you read a little about the 30 Years War, you
will almost certainly learn that in 1631 a certain large German
city was almost completely destroyed during that conflict. The
same German city has played a significant role in the development
of numerical methods for convection-diffusion problems. I refer to
{\em Magdeburg}.

Late 20th-century mathematicians from Magdeburg who have worked on
numerical methods for convection-diffusion problems include {\em
Goering,} Tobiska, Roos, Lube, Felgenhauer,
  John, Matthies, Risch, Schieweck, \ldots
Herbert Goering was the father of this school; all other names
here were his students or his students' students from Magdeburg.

\section{Where is our \lq\lq Peace of Westphalia"?}
Can we find a numerical method that is completely satisfactory for
all convection-diffusion problems? The general consensus seems to
be that in the future we will use adaptive meshes based on a
posteriori error indicators. Unfortunately the development of this
theory for convection-diffusion problems is only beginning; John
\cite{Jo00} numerically investigates several standard a posteriori
error indicators, and concludes that all are unsatisfactory to
varying degrees. Consequently I have not discussed this promising
line of attack in my talk.

To conclude, I would like to offer to young researchers embarking
on a study of convection-diffusion problems some advice drawn from
my own experience : always try the easiest case
first---it may be harder than you expect!

\end{document}